\documentclass{amsart}
\usepackage{amsmath}
\usepackage{amssymb}
\usepackage{amsthm}

\title{A modular branching rule for the generalized symmetric groups}

\author{Shunsuke Tsuchioka}
\address{Research Institute for Mathematical Sciences, %
         Kyoto University, Kyoto 606-8502. Japan.}
\email{tshun@kurims.kyoto-u.ac.jp}

\date{}

\newtheorem{Thm}{Theorem}[section]
\newtheorem{Def}[Thm]{Definition}
\newtheorem{DefThm}[Thm]{Definition and Theorem}

\newtheorem{SubLem}[Thm]{Sublemma}
\newtheorem{Rem}[Thm]{Remark}

\newcommand{\GE}{\mathfrak{g}}

\newcommand{\TE}{\widetilde{E}}
\newcommand{\TD}{\widetilde{D}}
\newcommand{\SEMI}{\wr}

\newcommand{\BA}{\overrightarrow{\boldsymbol{a}}}
\newcommand{\BB}{\overrightarrow{\boldsymbol{b}}}

\newcommand{\BC}{\overrightarrow{\boldsymbol{c}}}
\newcommand{\BN}{\overrightarrow{\boldsymbol{n}}}
\newcommand{\LAMBDA}{\overrightarrow{\boldsymbol{\lambda}}}
\newcommand{\MU}{\overrightarrow{\boldsymbol{\mu}}}
\newcommand{\NU}{\overrightarrow{\boldsymbol{\nu}_1}}
\newcommand{\NUU}{\overrightarrow{\boldsymbol{\nu}_2}}

\newcommand{\B}{\mathbb{B}}

\newcommand{\Z}{\mathbb{Z}}

\newcommand{\MH}{\mathcal{H}}

\newcommand{\II}{\mathfrak{I}}

\newcommand{\MP}{\mathcal{P}}

\newcommand{\SYM}[1]{\mathfrak{S}_{#1}}

\newcommand{\ISOM}{\stackrel{\sim}{\longrightarrow}}

\newcommand{\DEF}{\stackrel{\text{def}}{=}}

\DeclareMathOperator{\GL}{GL}

\DeclareMathOperator{\IRR}{\mathsf{Irr}}
\DeclareMathOperator{\HOM}{Hom}

\DeclareMathOperator{\AUT}{Aut}

\DeclareMathOperator{\CHAR}{char}

\DeclareMathOperator{\SOC}{\mathsf{Soc}}

\DeclareMathOperator{\RES}{\mathsf{Res}}

\DeclareMathOperator{\IND}{\mathsf{Ind}}

\begin{document}
\maketitle

\begin{abstract}
We give a modular branching rule for certain wreath products
as a generalization of Kleshchev's modular branching rule
for the symmetric groups.
Our result contains a modular branching rule for the complex reflection groups $G(m,1,n)$
(which are often called the generalized symmetric groups) in splitting fields for $\Z/m\Z$.
Especially for $m=2$ (which is the case of the Weyl groups of type $B$), we can give
a modular branching rule in any field.
Our proof is elementary in that it is essentially a combination of 
Frobenius reciprocity, 
Mackey theorem, 
Clifford's theory
and Kleshchev's modular branching rule.
\end{abstract}

\section{Introduction}
Given a sequence of groups such as
$
G_0\subseteq G_1\subseteq G_2\subseteq \cdots
$,
branching rule for this sequence is the rule that ``describes'' 
$\RES^{G_{n+1}}_{G_{n}}(V)$
for $G_{n+1}$-module $V$
or $\IND^{G_{n+1}}_{G_{n}}(W)$ 
for $G_{n}$-module $W$.
Let us review the case of the symmetric groups $\{G_n = \SYM{n}\}_{n\geq 0}$ with 
the default embedding of $\SYM{m}$ into $\SYM{n}$ for $m<n$, that is the embedding with respect
to the first $m$ letters.
In characteristic zero, we can summarize the classically known branching rule as follows
(for the details, see Theorem \ref{branching theorem}).
\begin{itemize}
\item for any irreducible $\SYM{n+1}$-module $V$, $\RES^{\SYM{n+1}}_{\SYM{n}}(V)$ is multiplicity-free.
\item Young's lattice controls the structure of $\RES^{\SYM{n+1}}_{\SYM{n}}(V)$ as $\SYM{n}$-module.
\end{itemize}
Here the meaning of the word ``multiplicity-free'' is not ambiguous because 
$\RES^{\SYM{n+1}}_{\SYM{n}}(V)$ is completely reducible in characteristic zero.

Recently, A.Kleshchev successfully discovered and proved its analogue in positive characteristics 
which is now known as Kleshchev's modular branching rule~\cite{Kl1, Kl2, Kl3, Kl4}.
The language of quantum groups and Kashiwara's crystal bases~\cite{HK,Kas}
lets us state Kleshchev's modular branching rule 
in chacateristic $p>0$ succinctly and beautifully as follows
(for the details, see Theorem \ref{modular branching theorem}).

\begin{itemize}
\item for any irreducible $\SYM{n+1}$-module $V$, $\SOC(\RES^{\SYM{n+1}}_{\SYM{n}}(V))$ 
is multiplicity-free.
\item the crystal basis $B(\Lambda_0)$ of the fundamental irreducible $U_q(\GE(A^{(1)}_{p-1}))$-module 
$L(\Lambda_0)$ controls the structure of $\SOC(\RES^{\SYM{n+1}}_{\SYM{n}}(V))$ as $\SYM{n}$-module.
\end{itemize}
Here for an $A$-module $X$ we denote by $\SOC(X)$ the largest completely reducible $A$-submodule of $X$.

We generalize the above Kleshchev's modular branching rule for the symmetric
groups to certain wreath products. 
Let $G$ be a finite group and $F$ be its splitting field and further assume 
that any irreducible $FG$-module is $1$-dimensional.
We denote by $\alpha$ the 
number of inequivalent irreducible representations of $G$.
For example, if $G$ is a $p$-group then $\alpha=1$.
The main result of this paper is succinctly stated as follows
(for the details, see Theorem \ref{main result}).

\begin{itemize}
\item for any irreducible $G\wr\SYM{n+1}$-module $V$, 
$\SOC(\RES^{G\wr\SYM{n+1}}_{G\wr\SYM{n}}(V))$ is multiplicity-free.
\item the crystal basis of 
$U_q(\GE(A^{(1)}_{p-1}))^{\otimes \alpha}(\cong U_q(\GE((A^{(1)}_{p-1})^{\oplus \alpha})))$-module 
$L(\Lambda_0)^{\otimes \alpha}$ controls the structure of 
$\SOC(\RES^{G\wr\SYM{n+1}}_{G\wr\SYM{n}}(V))$ as $G\wr\SYM{n}$-module.
\end{itemize}

Here for simplicity we assume $\CHAR F=p>0$, but our result also contains the case of
$\CHAR F=0$.
Note that our result contains a modular branching rule for the complex reflection groups 
$G(m,1,n)=(\Z/m\Z)\wr\SYM{n}$ 
(which are often called the generalized symmetric groups) in splitting fields for $\Z/m\Z$.
Especially for $m=2$ (which is the case of
the Weyl groups of type $B$), we can give
a modular branching rule in any field.

We recall a known generalization of Kleshchev's modular branching rule.
Let $\lambda$ be a positive integral weight of $U_q(\GE(A^{(1)}_{p-1}))$ and 
consider the inductive system of the cyclotomic degenerate Hecke algebras 
$\{\MH_n^{\lambda}\}_{n\geq 0}$~\cite[Chapter 7]{Kl6}.
Then the following holds~\cite{Kl6}. See also ~\cite{Ari,Bru,GV}.
If we take $\lambda=\Lambda_0$ then we get 
Kleshchev's modular branching rule for the symmetric groups.

\begin{itemize}
\item for any irreducible $\MH_{n+1}^{\lambda}$-module $V$, 
$\SOC(\RES^{\MH_{n+1}^{\lambda}}_{\MH_{n}^{\lambda}}(V))$ is multiplicity-free.
\item the crystal basis $B(\lambda)$ of 
$U_q(\GE(A^{(1)}_{p-1}))$-module 
$L(\lambda)$ controls the structure of 
$\SOC(\RES^{\MH_{n+1}^{\lambda}}_{\MH_{n}^{\lambda}}(V))$ as $\MH^{\lambda}_{n}$-module.
\end{itemize}

It should be noted that the different crystals appear in  
modular branching for the wreath products and for the Hecke algebras.
Although known proofs for the Hecke algebra case require the recent advances 
in modular representation 
theory actively involving other areas such as quantum groups, our case is elementary in that 
our proof is essentially a combination of Frobenius reciprocity, Mackey theorem, 
Clifford's theory and Kleshchev's modular branching rule.

\hspace{0mm}

\noindent\textbf{Acknowledgements}
The author would like to thank professor Masaki Kashiwara and 
Susumu Ariki for valuable comments.

\hspace{0mm}

\noindent\textbf{Notations and conventions}
In the following discussion of this paper, we assume for simplicity that 
any group is finite and any module is finite dimensional.

\begin{itemize}
\item For a finite dimensional algebra $A$, we denote by $\IRR(A)$ the set of 
isomorphism classes of irreducible $A$-modules.
\item For an $A$-module $V$, we denote by $[V]_A$ the isomorphism class of $A$-modules 
isomorphic to $V$. If $A$ is clear from the context, we often omit the suffix.
\item When we say directed graphs without modifiers, 
it means directed graphs with no loops nor 
no multiple-egdes. We write directed graph as $X=(V,E)$ where $V$ is the set of vertices
and $E\subseteq V\times V$ is its adjacent relation meaning there is a directed arrow
from $v_1$ to $v_2$ if and only if $(v_1,v_2)\in E$.
\item Let $l\geq 2$ be a positive integer.
A partition $(\lambda_1, \cdots, \lambda_k)$ is called $l$-restricted if
$\lambda_i-\lambda_{i+1}<l$ for any $1\leq i<k$. It is called $l$-regular if
its conjugate partition is $l$-restricted. 
All the partitions are defined to be 0-restricted and 0-regular.
\item Let $p>0$ be a prime number and $G$ be a group.
A conjugacy class of $G$ is called $p$-regular if 
the order of (one of or equivalently all of) its elements is prime to $p$.
All the congugacy classes of $G$ are defined to be 0-regular.
\item For a given sequence of non-negative integers 
$\BN = (n_1,\cdots,n_{\alpha})$ such that 
$\sum_{\beta=1}^{\alpha}n_{\beta}=n$ (such $\BN$ is called a composition of $n$), 
we denote by ${\SYM{}}_{\BN}$ the Young subgroup of $\SYM{n}$, that is
\begin{align*} 
{\SYM{}}_{\BN}
&\DEF
\SYM{\{1,\cdots,n_1\}}\times\SYM{\{n_1+1,\cdots,n_1+n_2\}}
\times\cdots\times\SYM{\{n-n_{\beta}+1,\cdots,n\}} \\
&\cong
\SYM{n_1}\times\cdots\times\SYM{n_{\beta}}.
\end{align*}
\end{itemize}

\section{Definitions}

\begin{Def}
Given a field $F$ and an inductive system of groups $\II=(\{G_n\}_{n\geq 0}, 
\{\varphi_n:G_n\to G_{n+1}\}_{n\geq 0})$, 
we define \textbf{modular branching graph} $\B_F(\II)$, which is 
a priori a directed graph with multiple-edges, as follows.
\begin{itemize}
\item the vertices are the elements of 
$\displaystyle\bigsqcup_{n\geq 0}\IRR(FG_n)$.
\item for two vertices 
$[W]\in\IRR(FG_n)$ and 
$[V]\in\IRR(FG_{n+1})$ of $\B_F(\II)$, there are
$\dim\HOM_{FG_n}(W, \RES^{FG_{n+1}}_{FG_n}V)$ edges from $[W]$ to $[V]$.
\end{itemize}
If $\B_F(\II)$ has only single-edges, 
we say that $\II$ is \textbf{socle multiplicity-free} over $F$.
\label{graph}
\end{Def}

\begin{Rem}
Note that for $FG_n$-module $W$ and $FG_{n+1}$-module $V$ we have
\begin{align*}
{} &{\quad}
\dim\HOM_{G_{n+1}}(\IND^{G_{n+1}}_{G_n}(W),V)
=
\dim\HOM_{G_n}(W,\RES^{G_{n+1}}_{G_n}(V)) \\
&=
\dim\HOM_{G_{n+1}}(V^{*}, \IND^{G_{n+1}}_{G_n}(W^{*}))
=
\dim\HOM_{G_n}(\RES^{G_{n+1}}_{G_n}(V^{*}), W^{*}) 
\end{align*}
by Frobenius reciprocity and Nakayama relations.
Hence the natural 4 choices for the definition of modular branching graph 
all coincide with each other if all the irreducible $FG_n$-modules are self-dual for each $n$.
\end{Rem}

\begin{Rem}
There is yet another natural definition of modular branching graph that
replaces $\dim\HOM_{G_n}(W,\RES^{G_{n+1}}_{G_n}V)$ by $[\RES^{G_{n+1}}_{G_n}V:W]$.
For the symmetric groups case $\{G_n = \SYM{n}\}_{n\geq 0}$, 
Kleshchev partially succeeded in describing this~\cite{Kl7} (and as a corollary we know that
the decomposition number of $\SYM{n}$ in positive characteristics can be arbitrary large), 
and Kleshchev also showed that knowing this type of modular branching graph is as hard as 
knowing the decomposition numbers~\cite{Kl5}.
\end{Rem}

\begin{Def}
Let $X_1=(V_1, E_1), \cdots, X_k=(V_k, E_k)$ be $k$ directed graphs.
We define the directed graph $X_1*\cdots*X_k=(V,E)$ as follows.
\begin{itemize}
\item $V = V_1\times \cdots\times V_k$.
\item $((v_1,\cdots,v_k),(w_1,\cdots,w_k))\in E$ iff there exists unique
$1\leq j\leq k$ such that $(v_j,w_j)\in E_j$ and $v_{j'}=w_{j'}$ for all $j'\ne j$.
\end{itemize}
Especially when $X_1=\cdots=X_k=X$, we write $X_1*\cdots*X_k=X^{*k}$.
\end{Def}

\section{Representation theory of the symmetric groups}

We shall first recall the representation theory of the symmetric groups.
It is well-known that for each partition $\lambda$ of $n$, we can construct
$\Z$-free, $\Z$-finite rank, 
$\Z\SYM{n}$-module $S^{\lambda}$ which is 
called the \textbf{Specht module}~\cite[Chapter 4]{Jam}.
Each $S^{\lambda}$ has an $\SYM{n}$-invariant symmetric bilinear form.
For any field $F$, we write $S^{\lambda}_F=F\otimes_\Z S^{\lambda}$ 
and denote by $D^{\lambda}_F$
the quotient of $S^{\lambda}_F$ by the radical of its invariant form.
It is also well-known 
that $D^{\lambda}_F=S^{\lambda}_F$ if $\CHAR F=0$~\cite[Chapter 4]{Jam}.
The following is the fundamental theorem of the representation theory of the 
symmetric groups.

\begin{Thm}[{\cite[Theorem 11.5]{Jam}}]
Suppose that our ground field $F$ has characteristic $p(\geq 0)$.
As $\mu$ varies over $p$-regular partitions of $n$, $D^{\mu}_F$ varies
over a complete set of inequivalent irreducible $F\SYM{n}$-modules.
Each $D^{\mu}_F$ is self-dual and absolutely irreducible. Every field is 
splitting field for $\SYM{n}$.
\label{symmetric}
\end{Thm}

\begin{Def}
We denote by $\SYM{}$ the inductive system of the symmetric groups
\begin{align*}
\SYM{} = (\SYM{0}\subseteq \SYM{1}\subseteq \SYM{2}\subseteq\cdots)
\end{align*}
Here if $m<n$, the default embedding of $\SYM{m}$ into $\SYM{n}$ is with respect
to the first $m$ letters.
\end{Def}

\begin{Thm}[{\cite[Theorem 9.2]{Jam}}]
For any field $F$ of characteristic zero, $\SYM{}$ is socle multiplicity-free over $F$.
Moreover, the map $\B_0\to\B_F(\SYM{}), \lambda\mapsto [S^{\lambda}]$ is an 
isomorphism as directed graphs where $\B_0=(\MP_0, E_0)$ is the directed graph obtained from 
the Young's lattice $\MP\DEF\{\lambda\vdash n\mid n\geq 0\}(=\MP_0)$ in the trivial manner.
\label{branching theorem}
\end{Thm}

Theorem \ref{branching theorem} is known as the classical branching rule.
Its analog in positive characteristics is known as Kleshchev's modular branching rule.
First we introduce the necessary terminology.

\begin{Def}
Let $l\geq 2$ be a positive integer and $i\in\Z/l\Z$.
Note that the following definitions depend on this given $l$.
Let $\lambda$ be a partition. 
We identify the partition $\lambda$ with the Young diagram of shape $\lambda$, 
that is $\{(r,s)\in\Z_{>0}\times\Z_{>0}|s\leq\lambda_r\}$.
\item[(1)] 
For a node $A=(r,s)\in\Z_{>0}\times\Z_{>0}$, 
we define its \textbf{residue} to be $-r+s+l\Z\in\Z/l\Z$.
\item[(2)]
A node $A$ inside $\lambda$ is called \textbf{$i$-removable} if
the residue of $A$ equals $i$ and $\lambda\setminus\{A\}$ is still a Young diagram.
\item[(3)]
A node $A$ outside $\lambda$ is called \textbf{$i$-addable} if
the residue of $A$ equals $i$ and $\lambda\cup\{A\}$ is again a Young diagram.
\item[(4)]
Now label all $i$-addable nodes of $\lambda$ by $+$ and all $i$-removable nodes of $\lambda$ by $-$.
The \textbf{$i$-signature} of $\lambda$ is the sequence of pluses and minuses obtained by
going along the rim of the Young diagram from bottom left to top right and reading off all
the signs.
\item[(5)] The \textbf{reduced $i$-signature} of $\lambda$ is the sequence of pluses and minuses 
obtained from the $i$-signature of $\lambda$ by successively erasing all neighboring pairs
of the form $-+$. Note that the reduced $i$-signature of $\lambda$ always looks like
a sequence that starts with $+$s followed by $-$s.
\item[(6)] Nodes that correspond to $+$s of 
the reduced $i$-signature of $\lambda$ are called \textbf{$i$-conormal}.
\item[(7)] The node that corresponds to the rightmost $+$ of 
the reduced $i$-signature of $\lambda$ is called \textbf{$i$-cogood}.
It is the rightmost $i$-conormal node.
\item[(8)] We set $\varphi_i(\lambda)=\#\{\textrm{$i$-conormal nodes of $\lambda$}\}$.
\item[(9)] If $\varphi_i(\lambda)>0$, we set $\tilde{f}_i(\lambda)=\lambda\cup\{A\}$
where $A$ is the (unique) $i$-cogood node.
\label{crystal terminology}
\end{Def}

\begin{Def}
Let $p$ be a prime number. We define the directed graph $\B_p = (\MP_p, E_p)$ where
$\MP_p = \{\lambda\in\MP\mid\textrm{$\lambda$ is $p$-regular}\}$ 
and the adjacent relation $E_p = 
\{(\lambda,\mu)\in\MP_p\times\MP_p\mid
\textrm{$\exists i\in\Z/p\Z, 
\varphi_i(\lambda)>0$ and 
$\tilde{f}_i(\lambda)=\mu$}\}$.
\end{Def}

The following is the form of Kleshchev's modular branching rule that we use.

\begin{Thm}[{\cite{Kl2,Kl6}}]
For any field $F$ whose characteristic is $p>0$, 
$\SYM{}$ is socle multiplicity-free over $F$.
Moreover, the map $\B_p\to\B_F(\SYM{}), \lambda\mapsto [D^{\lambda}_F]$ is an 
isomorphism as directed graphs.
\label{modular branching theorem}
\end{Thm}

\section{Representation theory of wreath products}

Let $G$ be a group and $F$ be a splitting field for $G$. 
We write $\IRR(FG)=\{[V_1],\cdots,[V_{\alpha}]\}$.
Definitions in this chapter depend on these datum.

Recall that the wreath product $G\wr\SYM{n}$ is the semi-direct product $G^n\rtimes_{\theta}\SYM{n}$
where 
\begin{align*}
\theta:\SYM{n}\longrightarrow \AUT(G^n),
\quad
\sigma\mapsto\theta(\sigma)((g_1,\cdots,g_n))=
(g_{\sigma^{-1}(1)},\cdots,g_{\sigma^{-1}(n)}).
\end{align*}
Hence any element $x\in G\wr\SYM{n}$ is written $x=(f;\pi)$ for uniquely determined
$f=(g_1,\cdots,g_n)\in G^n$ and $\sigma\in\SYM{n}$ and the multiplication rule is
\begin{align}
(f_1;\sigma_1)\cdot(f_2;\sigma_2) = (f_1\cdot(f_2)_{\sigma_1};\sigma_1\cdot\sigma_2)
\label{mult}
\end{align}
where $(f_2)_{\sigma_1} \DEF \theta(\sigma_1)(f_2)$.
The normal subgroup $(G^n\cong)\{(f;1_{\SYM{n}})\mid f\in G^n\}\subseteq G\wr\SYM{n}$ is often 
called the \textbf{base group} of $G\wr\SYM{n}$ and denoted by $G^*$.

The representation theory of wreath products is a typical application of Clifford's theory and
well-presented in \cite[Chapter 4]{JK}.
This chapter is a very brief summary of \cite[Chapter 4]{JK}.
Note that notations and the assumptions for definitions or theorems 
are slightly changed.

\begin{Def}
For a given composition $\BN = (n_1,\cdots,n_{\alpha})$ of $n$
and a permutation $\pi\in\SYM{\alpha}$, 
we define irreducible $FG^*$-module 
$E(\BN;\pi)=V_{\pi(1)}^{\otimes n_1}\otimes\cdots\otimes V_{\pi(\alpha)}^{\otimes n_{\alpha}}$.
We usually omit $\pi$ when $\pi = 1_{\SYM{\alpha}}$.
\end{Def}

\begin{Thm}[{\cite[4.3.27]{JK}}]
Let $\BN$ be as above.
The inertia group for $E(\BN)$ is given by $\{(f;\sigma)\mid f\in G^n, 
\sigma\in{\SYM{}}_{\BN}\}(\subseteq G\wr\SYM{n})$.
We denote this group by $G^*{\SYM{}}_{\BN}$.
\end{Thm}


\begin{Def}[{\cite[p.154]{JK}}]
Let $\BN$ and $\pi$ be as above.
To the underlying vector space $E(\BN;\pi)$, we can define a $F[G^*{\SYM{}}_{\BN}]$-module 
structure by
\begin{align*}
(f;\sigma)(v_1\otimes\cdots\otimes v_n)
\DEF
g_1v_{\sigma^{-1}(1)}\otimes\cdots\otimes g_nv_{\sigma^{-1}}(n)
\end{align*}
where $f = (g_1,\cdots,g_n)\in G^n$ and $\sigma\in {\SYM{}}_{\BN}$.
We denote this 
module 
by $\TE(\BN;\pi)$.
We usually omit $\pi$ when $\pi = 1_{\SYM{\alpha}}$.
\label{symmetric3}
\end{Def}

\begin{Def}
For a given sequence of $p$-regular 
partitions $\LAMBDA = (\lambda_1,\cdots,\lambda_{\alpha})$ such that 
$\sum_{\beta=1}^{\alpha}|\lambda_{\beta}|=n$.
We write $\BN=(|\lambda_1|,\cdots,|\lambda_{\alpha}|)$ and define a
irreducible $F{\SYM{}}_{\BN}$-module $D(\LAMBDA) = 
D^{\lambda_1}_F\otimes\cdots\otimes D^{\lambda_{\alpha}}_F$.
\end{Def}

\begin{Def}[{\cite[4.3.31]{JK}}]
Let $\LAMBDA$ and $\BN$ be as above.
To the underlying vector space $D(\LAMBDA)$, 
we define a $F[G^*{\SYM{}}_{\BN}]$-module 
structure by
\begin{align*}
(f;\sigma)(w_1\otimes\cdots\otimes w_{\alpha})
\DEF
\sigma(w_1\otimes\cdots\otimes w_{\alpha})
\end{align*}
where $f = (g_1,\cdots,g_n)\in G^n$ and $\sigma\in {\SYM{}}_{\BN}$.
We denote this 
module 
by $\TD(\LAMBDA)$.
\label{symmetric2}
\end{Def}

\begin{DefThm}[{\cite[4.4.3]{JK}}]
Let $\LAMBDA$ and $\BN$ be as above.
We define a $F[G\wr\SYM{n}]$-module $C(\LAMBDA)$ by
\begin{align*}
C(\LAMBDA) = \IND^{F[G\wr\SYM{n}]}_{F[G^*{\SYM{}}_{\BN}]}
\TE(\BN)\otimes \TD(\LAMBDA).
\end{align*}
If $\LAMBDA = (\lambda_1,\cdots,\lambda_{\alpha})\in\MP_p^{\alpha}$ varies while satisfying
$\sum_{\beta=1}^{\alpha}|\lambda_{\beta}|=n$, $C(\LAMBDA)$ varies
over a complete set of inequivalent irreducible $F[G\wr\SYM{n}]$-modules.
\end{DefThm}

\section{Main Results}
\begin{Def}
Denote by $G\wr\SYM{}$ the following inductive system of the wreath product groups of $G$.
\begin{align*}
G\wr\SYM{} = (G\wr\SYM{0}\subseteq G\wr\SYM{1}\subseteq G\wr\SYM{2}\subseteq\cdots).
\end{align*}
\end{Def}

\begin{Thm}
Let $G$ be a group and $F$ be its splitting filed of characteristic $p(\geq 0)$ such that
any irreducible $FG$-module is $1$-dimensional.
We denote by $\alpha$ the number 
of $p$-regular conjugacy classes of $G$ and fix a labeling of
$\IRR(FG)=\{[V_1],\cdots,[V_{\alpha}]\}$.
Then $G\SEMI\SYM{}$ is socle multiplicity-free over $F$.
Moreover, the map 
\begin{align*}
\B_p^{*\alpha}\longrightarrow\B_F(G\SEMI\SYM{}), 
\quad
(\lambda_1, \cdots, \lambda_{\alpha})\longmapsto [C(\lambda_1, \cdots, \lambda_{\alpha})]
\end{align*}
is an isomorphism as directed graphs.
\label{main result}
\end{Thm}

\begin{proof}
We show the equivalent statement that for any sequence of $p$-regular partitions
$\LAMBDA = (\lambda_1,\cdots, \lambda_{\alpha})$ and $\MU = (\mu_1,\cdots, \mu_{\alpha})$
such that $\sum_{i=1}^{\alpha}|\lambda_i|=n, \sum_{i=1}^{\alpha}|\mu_i|=n+1$, we have
\begin{align}
\dim\HOM_{G\wr\SYM{n}}(C(\LAMBDA),\RES^{G\wr\SYM{n+1}}_{G\wr\SYM{n}}C(\MU))\leq 1
\label{aim}
\end{align}
and the equality holds exactly when the following condition (\ref{cond}) holds.
\begin{align}
\textrm{there exists unique $1\leq\gamma\leq\alpha$ s.t.}
\begin{cases}
\textrm{$(\lambda_\gamma, \mu_\gamma)\in E_p$ (recall $\B_p = (\MP_p,E_p)$).} \\
\textrm{for any $\gamma'\ne\gamma$ we have $\lambda_{\gamma'}=\mu_{\gamma'}$.}
\end{cases}
\label{cond}
\end{align}

Let $\BA = (|\lambda_1|,\cdots,|\lambda_{\alpha}|), 
\BB = (|\mu_1|,\cdots,|\mu_{\alpha}|), 
\BC = (|\lambda_1|,\cdots,|\lambda_{\alpha}|, 1)$ and 
let $X = G\wr\SYM{n+1}, Y=G^{n+1}{\SYM{}}_{\BB}(\subseteq X), 
Z=G^{n}{\SYM{}}_{\BC}(\subseteq X), W = G\wr\SYM{n}$.
To avoid possible confusion, 
we reserve the trivial group isomorphism $t:(W\supseteq)G^{n}{\SYM{}}_{\BA}\ISOM Z$.
By Frobenius reciprocity, 
\begin{align*}
\dim\HOM_{W}(C(\LAMBDA),\RES^{X}_{W}C(\MU))
&=
\dim\HOM_{W}(\IND^{W}_{G^n{\SYM{}}_{\BA}}\TE(\BA)\otimes\TD(\LAMBDA),
\RES^{X}_{W}C(\MU)) \\
&=
\dim\HOM_{G^n{\SYM{}}_{\BA}}(\TE(\BA)\otimes\TD(\LAMBDA),
\RES^{W}_{G^n{\SYM{}}_{\BA}}\RES^{X}_{W}C(\MU)) \\
&=
\dim\HOM_{G^n{\SYM{}}_{\BA}}(\TE(\BA)\otimes\TD(\LAMBDA),
\RES^{X}_{G^n{\SYM{}}_{\BA}}C(\MU)).
\end{align*}

Let $D$ be a $(Z, Y)$-double coset representatives in $X$.
By Mackey theorem,
\begin{align*}
\RES^{X}_{G^n{\SYM{}}_{\BA}}C(\MU)
&=
{}^t\RES^{X}_{Z}C(\MU) \\
&=
{}^t\RES^{X}_{Z}\IND^{X}_{Y}(\TE(\BB)\otimes\TD(\MU)) \\
&\cong
\bigoplus_{d\in D} 
{}^t\IND^{Z}_{dYd^{-1}\cap Z}{}^{d}
(\RES^{Y}_{Y\cap d^{-1}Zd} \TE(\BB)\otimes\TD(\MU))
\end{align*}
as $F[G^n{\SYM{}}_{\BA}]$-modules 
where ${}^{d}M$ for $F[Y\cap d^{-1}Zd]$-module $M$ stands for a 
$F[dYd^{-1}\cap Z]$-module which is obtained by the pullback through the group isomorpshim
\begin{align*}
\varrho_d:dYd^{-1}\cap Z\ISOM Y\cap d^{-1}Zd,
\quad
x\mapsto d^{-1}xd
\end{align*}
and the same for $t$.

Now we recall a necessary fact about the double coset representatives of the symmetric groups.
Let $\NU,\NUU$ be compositions of $n$. Denote by $D_{\NUU}$ the set of minimal length
left ${\SYM{}}_{\NUU}$-coset representatives in $\SYM{n}$. Then $D_{\NU\NUU}\DEF
D_{\NU}^{-1}\cap D_{\NUU}$ is the set of minimal length 
$({\SYM{}}_{\NU},{\SYM{}}_{\NUU})$-double coset representatives in $\SYM{n}$~\cite{DJ}.
In the following discussion, it is important that we know $D_{\NUU}$ explicitly, that is
for $\NUU=(\eta_1,\cdots,\eta_{\kappa})$ we have
\begin{align}
D_{\NUU}=\left\{\sigma\in\SYM{n}\Bigg|
\substack{
\sigma(1)<\sigma(2)<\cdots<\sigma(\eta_1)
\\
\vdots\\
\sigma(n-\eta_{\kappa-1}+1)<\cdots<\sigma(n)}
\right\}.
\label{coset rep}
\end{align}
By the multiplication rule (\ref{mult}),
it is clear that 
$D_{\BC\BB}(\subseteq \SYM{n+1}\subseteq X)$ is a
$(Z,Y)$-double coset representatives in $X$.

So we need to compute for each $d\in D_{\BC\BB}$,
\begin{align*}
L(d)
\DEF
\dim\HOM_{G^n{\SYM{}}_{\BA}}\left(\TE(\BA)\otimes\TD(\LAMBDA),
{}^t\IND^{Z}_{dYd^{-1}\cap Z}{}^d(\RES^{Y}_{Y\cap d^{-1}Zd} \TE(\BB)\otimes\TD(\MU))
\right).
\end{align*}

We compute this value by taking their (various types of) duals.
Note that dual operation induces an involution $\IRR(FG)\to\IRR(FG)$.
We denote it by $\tau\in\SYM{\alpha}$ meaning 
that $[V_i^*]= [V_{\tau(i)}]$ for any $1\leq i\leq\alpha$.

\begin{SubLem}
\item[(P)]
$\TD(\LAMBDA)^{*}\cong\TD(\LAMBDA)$ as $F[G^n\SYM{\BA}]$-module.
\item[(Q)]
$\TE(\BA)^{*}\cong\TE(\BA;\tau)$ as $F[G^n\SYM{\BA}]$-module and
$\TE(\BB)^{*}\cong\TE(\BB;\tau)$ as $F[G^{n+1}\SYM{\BB}]$-module.
\end{SubLem}

\begin{proof}
\item[(P)]
By Theorem \ref{symmetric}, there exists $\SYM{|\lambda_k|}$-module isomorpshim 
$\Psi_{k}:D^{\lambda_k}_F\ISOM (D^{\lambda_k}_F)^{*}$ for each $1\leq k\leq\alpha$.
Then it is easy to check that the composition
\begin{align*}
D(\LAMBDA)
\xrightarrow[\Psi_1\otimes\cdots\otimes\Psi_{\alpha}]{\sim}
(D^{\lambda_1}_F)^{*}\otimes\cdots\otimes(D^{\lambda_{\alpha}}_F)^{*}
\xrightarrow[\textrm{can}]{\sim}
D(\LAMBDA)^{*}
\end{align*}
induces an $F[G^n\SYM{\BA}]$-module isomorphism
$\TD(\LAMBDA)\ISOM\TD(\LAMBDA)^{*}$
(see Definition \ref{symmetric2}).
\item[(Q)] The same as in (P) (see Definition \ref{symmetric3}).
\end{proof}

Therefore, for each $d\in D_{\BC\BB}$,
\begin{align}
\begin{split}
L(d)
&=
\dim\HOM_{G^n{\SYM{}}_{\BA}}\left(\TE(\BA)\otimes\TD(\LAMBDA),
{}^t\IND^{Z}_{dYd^{-1}\cap Z}{}^d(\RES^{Y}_{Y\cap d^{-1}Zd} \TE(\BB)\otimes\TD(\MU))
\right) \\
&=
\dim\HOM_{G^n{\SYM{}}_{\BA}}\left(
(
{}^t\IND^{Z}_{dYd^{-1}\cap Z}{}^d(
\RES^{Y}_{Y\cap d^{-1}Zd} \TE(\BB)\otimes\TD(\MU))
)^*,(\TE(\BA)\otimes\TD(\LAMBDA))^{*}
\right) \\
&=
\dim\HOM_{G^n{\SYM{}}_{\BA}}\left(
{}^t\IND^{Z}_{dYd^{-1}\cap Z}{}^d(
\RES^{Y}_{Y\cap d^{-1}Zd} \TE(\BB;\tau)\otimes\TD(\MU)),
\TE(\BA;\tau)\otimes\TD(\LAMBDA)
\right) \\
&=
\dim\HOM_{dYd^{-1}\cap Z}\left(
{}^d(
\RES^{Y}_{Y\cap d^{-1}Zd} \TE(\BB;\tau)\otimes\TD(\MU)),
\RES^{Z}_{dYd^{-1}\cap Z}{}^{t^{-1}}(\TE(\BA;\tau)\otimes\TD(\LAMBDA))
\right).
\end{split}
\label{final}
\end{align}

Note that $dYd^{-1}\cap Z=G^n(d{\SYM{}}_{\BB}d^{-1}\cap{\SYM{}}_{\BC})(\subseteq X)$.
By restricting to the subgroup $G^n\subseteq dYd^{-1}\cap Z$, we have
\begin{align*}
L(d)
&\leq
\dim\HOM_{G^n}\left(
\RES^{dYd^{-1}\cap Z}_{G^n}
{}^d(
\RES^{Y}_{Y\cap d^{-1}Zd} \TE(\BB;\tau)\otimes\TD(\MU)),
\RES^{Z}_{G^n}{}^{t^{-1}}(
\TE(\BA;\tau)\otimes\TD(\LAMBDA))
\right) \\
&=
\dim\HOM_{G^n}\left(
(V_{\tau(\xi(d^{-1}(1)))}\otimes\cdots\otimes V_{\tau(\xi(d^{-1}(n)))})^{\oplus\dim\TD(\MU)},
E(\BA;\tau)^{\oplus\dim\TD(\LAMBDA)}
\right)
\end{align*}
where we denote by $\xi(\chi)$ for $1\leq\chi\leq n+1$ the unique 
$1\leq\xi\leq\alpha$ such that
\begin{align*}
|\mu_1|+\cdots+|\mu_{\xi-1}|<\chi\leq|\mu_1|+\cdots+|\mu_{\xi}|.
\end{align*}

\begin{SubLem}
If $L(d)>0$, then there exists unique $1\leq j\leq\alpha$ such that
the followings are met.
\item[(A)] $|\lambda_{j'}|=|\mu_{j'}|$ for any $j'\ne j$.
\item[(B)] $|\lambda_j|+1=|\mu_j|$.
\item[(C)] $d(|\lambda_1|+\cdots+|\lambda_j|+1)=n+1$.
\item[(D)] $d^{-1}(1) < \cdots < d^{-1}(n)$.

Moreover, we have
\begin{align}
\begin{split}
U_{1,d}\DEF
d{\SYM{}}_{\BB}d^{-1}\cap {\SYM{}}_{\BC}
= 
{\SYM{}}_{(|\mu_1|,\cdots,|\mu_{j-1}|,|\mu_{j}|-1,|\mu_{j+1}|,\cdots,|\mu_{\alpha}|,1)}, \\
U_{2,d}\DEF
{\SYM{}}_{\BB}\cap d^{-1}{\SYM{}}_{\BC}d 
= 
{\SYM{}}_{(|\mu_1|,\cdots,|\mu_{j-1}|,|\mu_{j}|-1,1,|\mu_{j+1}|,\cdots,|\mu_{\alpha}|)}.
\label{perm3}
\end{split}
\end{align}
\label{perm2}
\end{SubLem}

\begin{proof}
We show that $j = \xi(d^{-1}(n+1))$.
Note that
\begin{align*}
{} &{}
\dim\HOM_{G^n}\left(
(V_{\tau(\xi(d^{-1}(1)))}\otimes\cdots\otimes V_{\tau(\xi(d^{-1}(n)))})^{\oplus\dim\TD(\MU)},
E(\BA;\tau)^{\oplus\dim\TD(\LAMBDA)}
\right)\\
&{=}
\dim\TD(\MU)\dim\TD(\LAMBDA)\dim\HOM_{G^n}\left(
(V_{\tau(\xi(d^{-1}(1)))}\otimes\cdots\otimes V_{\tau(\xi(d^{-1}(n)))}),
E(\BA;\tau)
\right)^{}\\
\end{align*}
and the isomorphim between $F$-vector spaces
\begin{align}
\begin{split}
{} &{}
\HOM_{G^n}\left(
(V_{\tau(\xi(d^{-1}(1)))}\otimes\cdots\otimes V_{\tau(\xi(d^{-1}(n)))}),
E(\BA;\tau)\right)\\
&{\cong}
\HOM_G(V_{\tau(\xi(d^{-1}(1)))}, V_{\tau(\zeta(1))})\otimes\cdots
\otimes\HOM_G(V_{\tau(\xi(d^{-1}(n)))}, V_{\tau(\zeta(n))}).
\end{split}
\label{isom1}
\end{align}
where we denote by $\zeta(\chi)$ for $1\leq\chi\leq n$ the unique 
$1\leq\zeta\leq\alpha$ such that
\begin{align*}
|\lambda_1|+\cdots+|\lambda_{\zeta-1}|<\chi\leq|\lambda_1|+\cdots+|\lambda_{\zeta}|.
\end{align*}

Hence, if $L(d)>0$, we have (by recalling $d\in D_{\BC}^{-1}\cap D_{\BB}\subseteq D_{\BC}^{-1}$)
\begin{align*}
\begin{cases}
1\leq 
{d^{-1}(1)<\cdots<d^{-1}(|\lambda_1|)} 
\leq |\mu_1| \\
|\mu_1|+1\leq 
{d^{-1}(|\lambda_1|+1)<\cdots<d^{-1}(|\lambda_1|+|\lambda_2|)}
\leq |\mu_1|+|\mu_2| \\
{\quad\quad\quad\quad\quad\quad\quad\quad\quad\quad\quad} {\vdots} {} \\
(n+1)-|\mu_\alpha|+1\leq 
{d^{-1}(n-|\lambda_\alpha|+1)<\cdots<d^{-1}(n)}
\leq n+1.
\end{cases}
\end{align*}
This implies that we have (A) and (B).
Hence it is enough to show that (C) holds.
Suppose to the contrary, we have
\begin{align*}
|\mu_1|+\cdots+|\mu_{j-1}|+1\leq 
d^{-1}(n+1)\leq |\mu_1|+\cdots+|\mu_j|-1.
\end{align*}
Because $d\in D_{\BC}^{-1}\cap D_{\BB}\subseteq D_{\BB}$, we have
\begin{align*}
d(|\mu_1|+\cdots+|\mu_{j-1}|+1)<
\cdots
<
d(d^{-1}(n+1))=n+1
<
\cdots
<d(|\mu_1|+\cdots+|\mu_j|).
\end{align*}
This is a contradiction.
Hence we have proved (A),(B),(C),(D)
(for unique $j = \xi(d^{-1}(n+1))$).
Since now we know the explicit form of $d$ characterized by (A),(B),(C),(D),
(\ref{perm3}) follows by the routine calculation.
\end{proof}

Now we assume $L(d)>0$ and $d$ be the form in Sublemma \ref{perm2} for 
uniquely determined $j$.
By restricting to the subgroup $U_{1,d}\subseteq dYd^{-1}\cap Z$, we have
\begin{align*}
L(d)
&\leq
\dim\HOM_{U_{1,d}}\left(
\RES^{dYd^{-1}\cap Z}_{U_{1,d}}
{}^d(
\RES^{Y}_{Y\cap d^{-1}Zd} \TE(\BB;\tau)\otimes\TD(\MU)),
\RES^{Z}_{U_{1,d}}
{}^{t^{-1}}(\TE(\BA;\tau)\otimes\TD(\LAMBDA))
\right) \\
&=
\dim\HOM_{U_{1,d}}\left(
{}^{\delta}(\RES^{{\SYM{}}_{\BB}}_{U_{2,d}}(D(\MU)^{\oplus\dim\TE(\BB;\tau)})),
\RES^{{\SYM{}}_{\BC}}_{U_{1,d}}({}^{T^{-1}}D(\LAMBDA)^{\oplus\dim\TE(\BA;\tau)})
\right) \\
&=
\dim\HOM_{U_{1,d}}\left(
{}^{\delta}(\RES^{{\SYM{}}_{\BB}}_{U_{2,d}}D(\MU)),
\RES^{{\SYM{}}_{\BC}}_{U_{1,d}}({}^{T^{-1}}D(\LAMBDA))
\right)
\end{align*}
where $T\DEF t|_{{\SYM{}}_{\BA}}:{\SYM{}}_{\BA}\ISOM{\SYM{}}_{\BC}$ and 
$\delta\DEF\varrho_{d}|_{U_{1,d}}:U_{1,d}\ISOM U_{2.d}$.
Note that here we use the fact that $\RES^{Y}_{U_{2,d}}(\TE(\BB;\tau))$ 
is a trivial $U_{2,d}$-module and 
$\RES^{Z}_{U_{1,d}}({}^{t^{-1}}\TE(\BA;\tau))$
is a trivial $U_{1,d}$-module because (any irreducible $FG$-module is $1$-dimensional and hence)
$\dim \TE(\BB;\tau)=\dim\TE(\BA;\tau)=1$
(see Definition \ref{symmetric2}).

By the explicit form of $d$ characterized by (A),(B),(C),(D) in Sublemma \ref{perm2},
and (\ref{perm3}), we have an isomorpshim as $F$-vector spaces
\begin{align}
\begin{split}
{} &{} 
\HOM_{U_{1,d}}\left(
{}^{\delta}(\RES^{{\SYM{}}_{\BB}}_{U_{2,d}}D(\MU)),
\RES^{{\SYM{}}_{\BC}}_{U_{1,d}}({}^{T^{-1}}D(\LAMBDA))
\right) \\
&\cong
\HOM_{\SYM{|\lambda_1|}}(D^{\mu_1}_F,D^{\lambda_1}_F)
\otimes\cdots\otimes
\HOM_{{\SYM{}}_{|\lambda_j|}}(
\RES^{\SYM{|\lambda_j|+1}}_{{\SYM{}}_{|\lambda_j|}}D^{\mu_j}_F,
D^{\lambda_j}_F) 
\otimes\cdots\otimes
\HOM_{\SYM{|\lambda_{\alpha}|}}(D^{\mu_{\alpha}}_F,D^{\lambda_{\alpha}}_F).
\end{split}
\label{isom2}
\end{align}

Applying (classical or Kleshchev's modular) branching rule for the symmetric groups, 
we have just proven that
\begin{itemize}
\item (\ref{aim}) holds.
\item if the equality for (\ref{aim}) holds then (\ref{cond}) holds.
\end{itemize}
So it remains to show that the equality for (\ref{aim}) holds if (\ref{cond}) holds.

\begin{SubLem}
Let $\mathcal{G}$ be a group and $\mathcal{G}_1, \mathcal{G}_2$ 
be its subgroups such that $\mathcal{G}_1\mathcal{G}_2=\mathcal{G}$ and 
$\mathcal{G}_1\cap \mathcal{G}_2=\{1_{\mathcal{G}}\}$.
Suppose we are given 4 representations of $G$
\begin{align*}
\begin{cases}
\rho_i:\mathcal{G}\longrightarrow\GL_F(\mathcal{V}_i) & (i=1,2) \\
\psi_i:\mathcal{G}\longrightarrow\GL_F(\mathcal{W}_i) & (i=1,2) \\
\end{cases}
\end{align*}
such that for any $g_1\in \mathcal{G}_1$ and $g_2\in \mathcal{G}_2$, we have
\begin{align*}
\begin{cases}
\rho_1(g_1g_2)=\rho_1(g_1),
&
\rho_2(g_1g_2)=\rho_2(g_2) \\
\psi_1(g_1g_2)=\psi_1(g_1),
&
\psi_2(g_1g_2)=\psi_2(g_2).
\end{cases}
\end{align*}
Then we have the inequality
\begin{align*}
{} &{}
\dim\HOM_{\mathcal{G}}(\mathcal{V}_1\otimes \mathcal{V}_2, 
\mathcal{W}_1\otimes \mathcal{W}_2) \\
&\geq
\dim\HOM_{\mathcal{G}_1}(\RES^{\mathcal{G}}_{\mathcal{G}_1}\mathcal{V}_1, 
\RES^{\mathcal{G}}_{\mathcal{G}_1}\mathcal{W}_1)
\cdot
\dim\HOM_{\mathcal{G}_2}(\RES^{\mathcal{G}}_{\mathcal{G}_2}\mathcal{V}_2, 
\RES^{\mathcal{G}}_{\mathcal{G}_2}\mathcal{W}_2).
\end{align*}
\label{geq}
\end{SubLem}

\begin{proof}
Note that there is a natural injection (between $F$-vector spaces)
\begin{align*}
{} &{}
\HOM_{\mathcal{G}_1}(\RES^{\mathcal{G}}_{\mathcal{G}_1}\mathcal{V}_1, 
\RES^{\mathcal{G}}_{\mathcal{G}_1}\mathcal{W}_1)
\otimes
\HOM_{\mathcal{G}_2}(\RES^{\mathcal{G}}_{\mathcal{G}_2}\mathcal{V}_2, 
\RES^{\mathcal{G}}_{\mathcal{G}_2}\mathcal{W}_2) \\
&\hookrightarrow
\HOM_{\mathcal{G}}(\mathcal{V}_1\otimes \mathcal{V}_2, 
\mathcal{W}_1\otimes \mathcal{W}_2)
\end{align*}
that sends $\varphi_1\otimes\varphi_2$ 
to $\varphi_1\otimes\varphi_2$.
\end{proof}

Let us assume that (\ref{cond}) holds.
Put $j=\gamma$ and take $d\in D_{\BC,\BB}$ characterized by (A),(B),(C),(D) in Sublemma \ref{perm2}.
As in the above discussion, we have only to show that $L(d)=1$.
Apply Sublemma \ref{geq} under 
\begin{align*}
\begin{cases}
\mathcal{G}=dYd^{-1}\cap Z=G^nU_{1,d}(\subseteq X),\quad
\mathcal{G}_1=G^n(\subseteq\mathcal{G}),
\quad
\mathcal{G}_2=U_{1,d}(\subseteq\mathcal{G}) \\
\mathcal{V}_1 = {}^{d}(\RES^{Y}_{d^{-1}\mathcal{G}d}\TE(\BB;\tau)),\quad
\mathcal{V}_2 = {}^{d}(\RES^{Y}_{d^{-1}\mathcal{G}d}\TD(\MU)) \\
\mathcal{W}_1 = \RES^{Z}_{\mathcal{G}}({}^{t^{-1}}\TE(\BA;\tau)),\quad
\mathcal{W}_2 = \RES^{Z}_{\mathcal{G}}({}^{t^{-1}}\TD(\LAMBDA))) \\
\end{cases}
\end{align*}
and we have
\begin{align*}
L(d)
&{\stackrel{\textrm{(\ref{final})}}{=}}
\dim\HOM_{dYd^{-1}\cap Z}\left(
{}^d(
\RES^{Y}_{Y\cap d^{-1}Zd} \TE(\BB;\tau)\otimes\TD(\MU)),
\RES^{Z}_{dYd^{-1}\cap Z}{}^{t^{-1}}(\TE(\BA;\tau)\otimes\TD(\LAMBDA))
\right) \\
&= \dim\HOM_{\mathcal{G}}(\mathcal{V}_1\otimes \mathcal{V}_2, 
\mathcal{W}_1\otimes \mathcal{W}_2) \\
&\geq
\dim\HOM_{\mathcal{G}_1}(\RES^{\mathcal{G}}_{\mathcal{G}_1}\mathcal{V}_1, 
\RES^{\mathcal{G}}_{\mathcal{G}_1}\mathcal{W}_1)
\cdot
\dim\HOM_{\mathcal{G}_2}(\RES^{\mathcal{G}}_{\mathcal{G}_2}\mathcal{V}_2, 
\RES^{\mathcal{G}}_{\mathcal{G}_2}\mathcal{W}_2) 
{\stackrel{\textrm{(\ref{isom1}),(\ref{isom2})}}{=}} 1.
\end{align*}

Because 
the converse inequality (\ref{aim}) has been already established,
we reach to the conclusion that $L(d)=1$.
\end{proof}

\end{document}